\documentclass{article}
\usepackage[utf8]{inputenc}
\usepackage{IEEEtrantools}
\usepackage{amsmath}
\usepackage{amssymb}
\usepackage{hyperref}
\usepackage{graphicx}
\usepackage[acronym,nomain]{glossaries}
\usepackage{booktabs} 
\usepackage{geometry} 
\usepackage{array} 
\newcolumntype{L}{>{\displaystyle}l}

\usepackage{fancyhdr} 
\usepackage{lscape}
\fancypagestyle{mylandscape}{
\fancyhf{} 
\fancyfoot{
\makebox[\textwidth][r]{
  \rlap{\hspace{.75cm}
    \smash{
      \raisebox{4.87in}{
        \rotatebox{90}{\thepage}}}}}}
}

\newacronym{BFS}{BFS}{Backward-forward sweep}
\newacronym{CIM}{CIM}{Current injection method}
\newacronym{PMD}{PMD}{PowerModelsDistribution.jl}
\newacronym{EN}{EN}{Explicity neutral}

\DeclareMathOperator{\diagm}{diagm}

\title{On the Implementation of the Fixed Point Iteration Current Injection Method to Solve Four-Wire Unbalanced Power Flow in \textsc{PowerModelsDistribution.jl}}

\author{Frederik Geth \thanks{\href{mailto:frederik.geth@gridqube.com}{frederik.geth@gridqube.com}. F. Geth carried out part of this work at CSIRO but has since joined GridQube.}
\and Sander Claeys \thanks{S. Claeys carried out part of this work while holding a doctoral (PhD) strategic basic research grant (1S82518N) of Research
Foundation - Flanders (FWO).}
\and Rahmat Heidari \thanks{\href{mailto:rahmat.heidarihaei@csiro.au}{rahmat.heidarihaei@csiro.au}}}

\date{February 2023}

\begin{document}

\maketitle

\section*{Abstract}
This report serves as a technology description of a Julia-based re-implementation of the fixed-point current injection algorithm, available in \textsc{PowerModelsDistribution.jl} \cite{FOBES2020106664}.
This report does \emph{not} describe a \emph{novel} method for solving unbalanced power flow problems.
It merely provides a description of the fixed point iteration variant of the current injection method, inspired  by the existing open-source implementation in \textsc{OpenDSS}\footnote{available:  \url{https://www.epri.com/pages/sa/opendss?lang=en}} \cite{Dugan6039829}. 
The current injection method is commonly conceived as a system of nonlinear equalities solved by Newton's method \cite{Araujo2006_edit,Garcia2000}.
However, as Roger Dugan points out in the \textsc{OpenDSS} documentation, the fixed point iteration variant commonly outperforms most  methods, while supporting meshed topologies from the ground up. 

We note that the unbalanced power flow algorithm in turn relies on matrix solvers for sparse systems of equations. 
In the context of circuits and factorizing nodal admittance matrices, the sparsity-exploiting `KLU' solver \cite{PalamadaiNatarajan2005} has proven to be both reliable and scalable. 
\textsc{OpenDSS}  uses KLU. 

This report documents work-in-progress, and the authors aim to update it when  learnings are obtained or more features are added to the implementation in \textsc{PowerModelsDistribution.jl}. 
The authors invite collaborators to contribute through pull requests on the repository\footnote{\url{https://github.com/lanl-ansi/PowerModelsDistribution.jl/pulls}}.

\newpage
\tableofcontents

\newpage

\section{Introduction on Unbalanced Power Flow}
Unbalanced power flow is the problem of finding a solution to the steady-state physics of ac power networks. Specifically, unbalanced power flow engines embed the multiconductor ac circuit laws, where series impedances and shunt admittances are matrices, not scalars.

Without aiming to do a detailed literature review,
\emph{unbalanced power flow algorithms} commonly discussed in the literature are:
 \begin{itemize}
     \item \gls{BFS} method, which is also a fixed-point iteration algorithm \cite{Cheng387902}.
     \item \gls{CIM} based on Newton-Raphson, using derivatives \cite{Penido1397701}.
     \item Current injection method as a fixed point iteration algorithm \cite{Dugan6039829}\footnote{OpenDSS setting \texttt{Algorithm=Normal}. Note that the Gauss-Seidel power flow method, often used in teaching power flow solution methods for transmission networks, is also a fixed-point iteraiton method, however distinct from OpenDSS's approach.}.
     \item Current injection method Newton-style as implemented in OpenDSS (also fixed-point, derivative-free) \cite{Dugan6039829} \footnote{OpenDSS setting \texttt{Algorithm=Newton}}.
 \end{itemize}
 
 Note that current injection type approaches are fundamentally centered around a nodal admittance representation of the circuit.
 Conversely, the canonical \gls{BFS} approach uses element-wise impedance-based representations and requires a radial topology\footnote{i.e. \gls{BFS} depends on hacks to overcome this limitation.}. 
The \gls{BFS} and \gls{CIM} fixed point algorithms have the following similarities:
 \begin{itemize}
     \item both are fixed-point iteration methods, and derivative-free,
     \item both are  flexible in the kind of load models they can accept, i.e.  power as a function of voltage does not need to be smooth.
 \end{itemize}
 
 Some of their key differences are:
 \begin{itemize}
     \item \gls{BFS} does not require any large matrix  factorization/solves, 
     \item fixed-point \gls{CIM} requires one factorization overall, and one matrix solve in each iteration,
     \item \gls{CIM} can handle meshed networks naturally, whereas \gls{BFS} requires workarounds.
 \end{itemize}
 Other algorithms  have  been proposed:
 \begin{itemize}
     \item Holomorphic embedding power flow for three-phase networks \cite{su11061774}.
     \item Interior point methods applied to a system of nonlinear equations representing the power flow physics \cite{FOBES2020106664}.
 \end{itemize}

 \subsection{Document Structure}
The document is structured as follows:
 \begin{itemize}
     \item \S2 introduces the notation and describes at an abstract level how components can be modeled.
     \item \S3 describes the fixed-point current injection algorithm, building on the mathematical models of the components, based on a system nodal admittance representation.
     \item \S4 derives the mathematical models for specific components that are supported, e.g. loads, generators, lines, cables, transformers.
     \item \S5 discusses the validation of the algorithm implementation w.r.t. OpenDSS.
 \end{itemize}

\section{Abstract Component Models}
Every network component $c$ -- including both power delivery and power consumption/injection --  connected to a set of bus-terminals pairs is modeled as a parallel composition of an admittance $\mathbf{Y}_c $ and a current source $\mathbf{I}^{\text{nl}}_c$,
  \begin{IEEEeqnarray}{C} 
\underbrace{\mathbf{I}_c}_{\text{component current }}  = \underbrace{\mathbf{Y}_c\mathbf{U}_c}_{\text{linear part}}  + \underbrace{ \mathbf{I}^{\text{nl}}_c. }_{\text{correction}} \label{eq_nonlinear_framework}
\end{IEEEeqnarray}

In this report, we illustrate expressions for a four-terminal/four-conductor network, which is general enough to capture up to four-wire networks with explicit neutral, as well as Kron-reduced four-wire networks. 
Nevertheless, the abstractions work for any number of terminals.
For instance, the NEVTestCase set up in \textsc{OpenDSS} contains up to 17 electromagnetically coupled wires.
The current and voltage vector sizes for different elements in the four-wire network are:
\begin{itemize}
    \item as power consumption/injection devices connect to a specific bus, the vectors $\mathbf{U}_c, \mathbf{I}_c$ have length 4,
    \item for power delivery elements, which connect two buses, the vectors $\mathbf{U}_c, \mathbf{I}_c$ have length $2\times 4 = 8$.
\end{itemize}

Note that the vectors do not need to be standardized to length 4. 
Nevertheless, in the description of models we assume this, to simplify notation. 
In the implementation, we use variable-size primitive vectors and matrices.

\subsection{Indices and Sets}
The followings sets and indices are used throughout the report:
  \begin{IEEEeqnarray}{L?s} 
i \in \mathcal{I} & buses \\
p \in \mathcal{P} \subseteq \{a,b,c,n\} & terminals \\
c \in \mathcal{C} & components \\
l \in \mathcal{L} \subset \mathcal{C} & lines \\
x \in \mathcal{X} \subset \mathcal{C} & transformers \\
(i,p) \in \mathcal{T}^{\text{bt}} \subseteq \mathcal{I} \times \mathcal{P} & bus-terminal topology\\
(c,i) \in \mathcal{T}^{\text{bus}} \subseteq \mathcal{C} \times \mathcal{I} & component-bus topology\\
(c,i,p) \in \mathcal{T}^{\text{term}} \subseteq \mathcal{T}^{\text{bus}} \times \mathcal{P} & component-bus-terminal topology
\end{IEEEeqnarray}

\subsection{Buses and Terminals}
A bus $i$ has terminals $p \in \{a,b,c,n\}$. 
The voltage phasor of bus $i$ assembles all the terminal voltages in the component-bus topology as,
\begin{IEEEeqnarray}{C} 
\underbrace{\mathbf{U}_c}_{\text{voltages seen by component $c$}} \leftarrow
\underbrace{\mathbf{U}_i}_{\text{voltages on bus $i$}} = 
\begin{bmatrix}
U_{i,a} \\
U_{i,b} \\
U_{i,c} \\
U_{i,n} \\
\end{bmatrix} \in \mathbb{C}^{4\times 1}, 
\quad (c,i) \in \mathcal{T}^{\text{bus}}.
\end{IEEEeqnarray}

Two types of reference nodes occur in the current version, which have a fixed voltage associated with them:
\begin{itemize}
    \item reference buses, where the phasor is fixed (all phase terminals of a reference bus are given known values)
    \item the neutral terminal of any bus when it is perfectly grounded (only the neutral of the bus is given a 0 value)
\end{itemize}

\subsection{Power Delivery Elements} \label{sec_power_del_framework}
\begin{figure}[tbh]
  \centering
    \includegraphics[width=0.90\columnwidth]{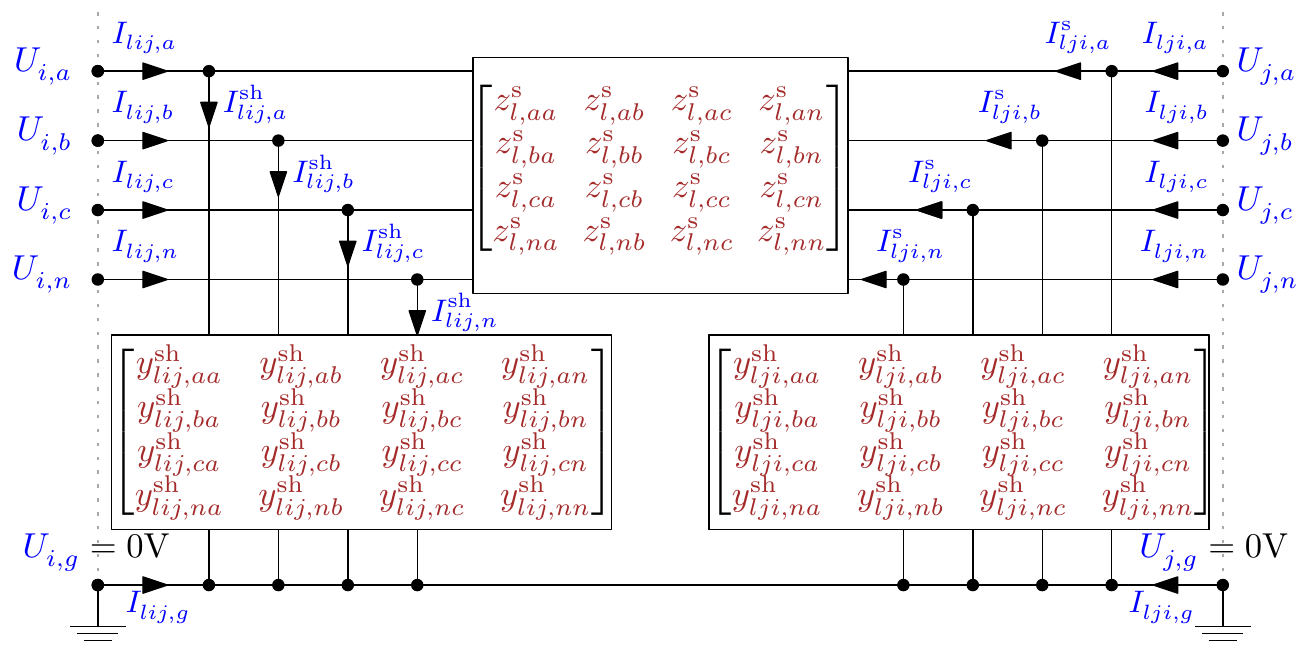}
  \caption{Four-wire branches have 8 nodes}  \label{fig_pd_element_4w}
\end{figure}

A 4-wire power delivery element such as branch $l$ connecting buses $i$ and $j$, shown in Figure~\ref{fig_pd_element_4w}, can be represented by the component-bus-terminal topology
\begin{IEEEeqnarray}{C} 
\{ (l,i,a), (l,i,b), (l,i,c), (l,i,n), (l,j,a), (l,j,b), (l,j,c), (l,j,n) \} \in \mathcal{T}^{\text{term}}.
\end{IEEEeqnarray}

This power delivery element has self and mutual impedances between all of its conductors. 
We define the matrix series admittance $ \mathbf{Y}^{\text{s}}_l$ of branch $l$, as
  \begin{IEEEeqnarray}{C} 
  \mathbf{Y}^{\text{s}}_l \leftarrow 
  \begin{bmatrix}
  y^{\text{s}}_{l,aa} &   y^{\text{s}}_{l,ab} &   y^{\text{s}}_{l,ac} &   y^{\text{s}}_{l,an} \\
  y^{\text{s}}_{l,ba} &   y^{\text{s}}_{l,bb} &   y^{\text{s}}_{l,bc} &   y^{\text{s}}_{l,bn} \\
  y^{\text{s}}_{l,ca} &   y^{\text{s}}_{l,cb} &   y^{\text{s}}_{l,cc} &   y^{\text{s}}_{l,cn} \\
  y^{\text{s}}_{l,na} &   y^{\text{s}}_{l,nb} &   y^{\text{s}}_{l,nc} &   y^{\text{s}}_{l,nn} \\
  \end{bmatrix} \in \mathbb{C}^{4\times 4}
  \end{IEEEeqnarray}
The sending end shunt admittance matrix $ \mathbf{Y}^{\text{sh}}_{lij}$ is
   \begin{IEEEeqnarray}{C} 
   \mathbf{Y}^{\text{sh}}_{lij} \leftarrow  
  \begin{bmatrix}
  y^{\text{sh}}_{lij,aa} &   y^{\text{sh}}_{lij,ab} &   y^{\text{sh}}_{lij,ac} &   y^{\text{sh}}_{lij,an} \\
  y^{\text{sh}}_{lij,ba} &   y^{\text{sh}}_{lij,bb} &   y^{\text{sh}}_{lij,bc} &   y^{\text{sh}}_{lij,bn} \\
  y^{\text{sh}}_{lij,ca} &   y^{\text{sh}}_{lij,cb} &   y^{\text{sh}}_{lij,cc} &   y^{\text{sh}}_{lij,cn} \\
  y^{\text{sh}}_{lij,na} &   y^{\text{sh}}_{lij,nb} &   y^{\text{sh}}_{lij,nc} &   y^{\text{sh}}_{lij,nn} \\
  \end{bmatrix} \in \mathbb{C}^{4\times 4},
    \end{IEEEeqnarray}
    and the receiving end shunt admittance matrix $ \mathbf{Y}^{\text{sh}}_{lji} $ is
   \begin{IEEEeqnarray}{C} 
   \mathbf{Y}^{\text{sh}}_{lji} \leftarrow  
  \begin{bmatrix}
  y^{\text{sh}}_{lji,aa} &   y^{\text{sh}}_{lji,ab} &   y^{\text{sh}}_{lji,ac} &   y^{\text{sh}}_{lji,an} \\
  y^{\text{sh}}_{lji,ba} &   y^{\text{sh}}_{lji,bb} &   y^{\text{sh}}_{lji,bc} &   y^{\text{sh}}_{lji,bn} \\
  y^{\text{sh}}_{lji,ca} &   y^{\text{sh}}_{lji,cb} &   y^{\text{sh}}_{lji,cc} &   y^{\text{sh}}_{lji,cn} \\
  y^{\text{sh}}_{lji,na} &   y^{\text{sh}}_{lji,nb} &   y^{\text{sh}}_{lji,nc} &   y^{\text{sh}}_{lji,nn} \\
  \end{bmatrix} \in \mathbb{C}^{4\times 4}.
    \end{IEEEeqnarray}
We can now construct the primitive branch admittance matrix $\mathbf{Y}^{\text{tot}}_l $ as
   \begin{IEEEeqnarray}{C} 
     \label{eq_Y_lines}
\mathbf{Y}^{\text{tot}}_l \leftarrow  
\begin{bmatrix}
\mathbf{Y}^s_l + \mathbf{Y}_{lij}^{\text{sh}} & - \mathbf{Y}^s_l \\
-\mathbf{Y}^s_l                  & \mathbf{Y}^s_l + \mathbf{Y}_{lji}^{\text{sh}}
\end{bmatrix} \in \mathbb{C}^{8\times 8}.
\end{IEEEeqnarray}
Note that on the basis of $\mathbf{Y}^s_l $ being rank-4, the rank of the whole matrix depends on that of the shunt terms. 
If $\mathbf{Y}_{lij}^{\text{sh}} = \mathbf{Y}_{lij}^{\text{sh}} = \mathbf{0}$, the resulting matrix has rank 4, but otherwise rank 8.

We compose a combined voltage vector $\mathbf{U}_l $ by stacking the sending $\mathbf{U}_i$ and receiving end $\mathbf{U}_j$ voltage phasors, 
  \begin{IEEEeqnarray}{C} 
\mathbf{U}_l = 
\begin{bmatrix}
\mathbf{U}_i\\
\mathbf{U}_j
\end{bmatrix}
=
\begin{bmatrix}
U_{i,a} \\
U_{i,b} \\
U_{i,c} \\
U_{i,n} \\
U_{j,a} \\
U_{j,b} \\
U_{j,c} \\
U_{j,n} \\
\end{bmatrix} \in \mathbb{C}^{8\times 1}.
\end{IEEEeqnarray}
Similarly, we compose the branch current vector $\mathbf{I}_l$ by stacking the sending $\mathbf{I}_{lij}$ and receiving $\mathbf{I}_{lji}$ current phasors, 
  \begin{IEEEeqnarray}{C} 
\mathbf{I}_l = 
\begin{bmatrix}
\mathbf{I}_{lij}\\
\mathbf{I}_{lji}
\end{bmatrix}
=
\begin{bmatrix}
I_{lij,a} \\
I_{lij,b} \\
I_{lij,c} \\
I_{lij,n} \\
I_{lji,a} \\
I_{lji,b} \\
I_{lji,c} \\
I_{lji,n} 
\end{bmatrix} \in \mathbb{C}^{8\times 1}.
\end{IEEEeqnarray}
Finally, the bus injection model is stated in matrix form as,
  \begin{IEEEeqnarray}{C} 
\mathbf{I}_l = 
\begin{bmatrix}
\mathbf{I}_{lij}\\
\mathbf{I}_{lji}
\end{bmatrix}
=
\begin{bmatrix}
\mathbf{Y}^{\text{s}}_l + \mathbf{Y}_{lij}^{\text{sh}} & - \mathbf{Y}^{\text{s}}_l \\
-\mathbf{Y}^{\text{s}}_l & \mathbf{Y}^{\text{s}}_l + \mathbf{Y}_{lji}^{\text{sh}}
\end{bmatrix}
\begin{bmatrix}
\mathbf{U}_i\\
\mathbf{U}_j
\end{bmatrix}
=
\mathbf{Y}^{\text{tot}}_l
\mathbf{U}_l.
\end{IEEEeqnarray}

For the circuit components that are linear in current-voltage variables, the compensation current is set to 0,
  \begin{IEEEeqnarray}{C} 
    \label{eq_Ic_lines}
\mathbf{I}^{\text{nl}}_c = \mathbf{0}.
\end{IEEEeqnarray}

\subsection{Power Consumption/Injection Elements}
At a 4-terminal bus, a load or generator (or storage system, ...) connected to bus $i$  can be represented by the component-bus-terminal topology
\begin{IEEEeqnarray}{C} 
\{ (c,i,a), (c,i,b), (c,i,c), (c,i,n) \} \in \mathcal{T}^{\text{term}}.
\end{IEEEeqnarray} 

\begin{figure}[tbh]
  \centering
    \includegraphics[width=0.5\columnwidth]{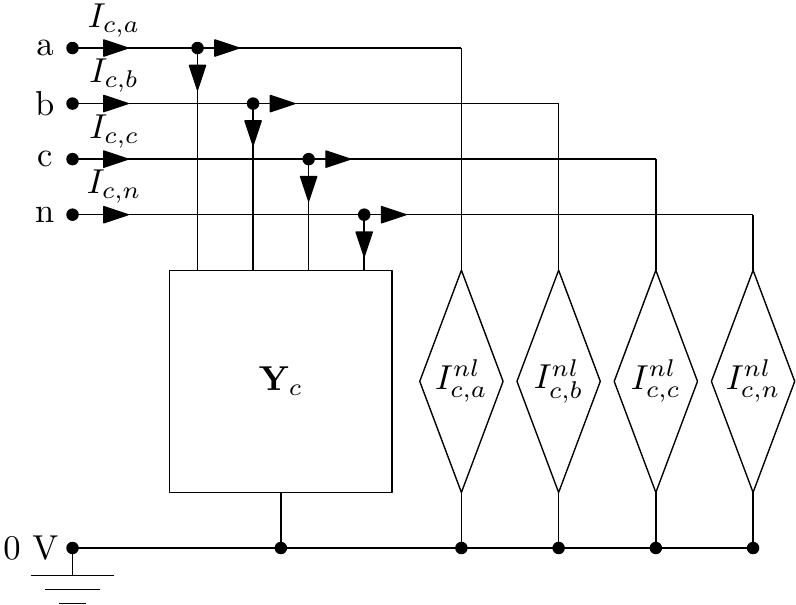}
  \caption{Nonlinear elements use a linear element in parallel with a controllable current source}  \label{fig_pd_element}
\end{figure}

The component is first approximated by a primitive admittance matrix $\mathbf{Y}_c$ (see Fig. \ref{fig_pd_element}), with a size proportional to the number of terminals it connects to:
  \begin{IEEEeqnarray}{C} 
  \mathbf{Y}_c \in \mathbb{C}^{4\times 4}.
\end{IEEEeqnarray}

Power consumption and generation devices, such as the ``constant power'' load model, often have nonlinear characteristics in the current-voltage space. 
Therefore, we iteratively adapt the ``compensation current'' injections $\mathbf{I}^{\text{nl}}_c$ to obtain the desired solution and solve equation \eqref{eq_nonlinear_framework}.
  \begin{IEEEeqnarray}{C} 
\mathbf{I}^{\text{nl}}_c = 
\begin{bmatrix}
I^{\text{nl}}_{c,a} \\
I^{\text{nl}}_{c,b} \\
I^{\text{nl}}_{c,c} \\
I^{\text{nl}}_{c,n} \\
\end{bmatrix} = f(\mathbf{U}_i) \in \mathbb{C}^{4\times 1}.
\end{IEEEeqnarray}

\section{Fixed Point Iteration Algorithm}

\subsection{Abstracting Network Details}

After building a model for each component in the network with the following structure:
  \begin{IEEEeqnarray}{C} 
\mathbf{I}_c  =\mathbf{Y}_c\mathbf{U}_c  + \mathbf{I}^{\text{nl}}_c, 
\qquad 
\mathbf{Y_c}: (i,p)\in \mathcal{T}^{\text{bt}} \rightarrow (c,i,p) \in \mathcal{T}^{\text{term}}
\qquad
\text{for }  (c,i)\in \mathcal{T}^{\text{bus}}.
\label{eq_nonlinear_framework_repeat}
\end{IEEEeqnarray}
we construct the \emph{system} nodal admittance matrix $\mathbf{Y}: \mathcal{T}^{\text{bt}} \rightarrow \mathcal{T}^{\text{term}}$ based on the primitive admittance matrices $\mathbf{Y}_c$ of the individual components $c$, by starting from a zero matrix, and adding the values of the primitive admittances to the appropriate rows and columns in the system matrix corresponding to the topological sets.

We partition the \emph{system} voltage vector $\mathbf{U}$ into \textbf{f}ixed (reference bus, grounded neutral terminals) and \textbf{v}ariable voltages, and derive the corresponding partitions for impedance and \emph{system} current variables in $\mathbf{I}$,
  \begin{IEEEeqnarray}{C} 
  \underbrace{
 \begin{bmatrix}
 \mathbf{I}^{\text{f}}\\
\mathbf{I}^{\text{v}}
 \end{bmatrix}
 }_{\mathbf{I}}
 =
   \underbrace{
 \begin{bmatrix}
  \mathbf{Y}^{\text{ff}}  & \mathbf{Y}^{\text{fv}} \\
  \mathbf{Y}^{\text{vf}}  & \mathbf{Y}^{\text{vv}}
 \end{bmatrix} 
  }_{\mathbf{Y}}
     \underbrace{
  \begin{bmatrix}
 \mathbf{U}^{\text{f}}\\
\mathbf{U}^{\text{v}}
 \end{bmatrix}
   }_{\mathbf{U}},
   \qquad
   \mathbf{Y}: \mathcal{T}^{\text{bt}} \rightarrow \mathcal{T}^{\text{term}}.
\end{IEEEeqnarray}
We define the variable current vector $\mathbf{I}^{\text{v}}$ as a vector of the nonlinear correction currents, $ \mathbf{I}^{\text{nl}}$, which is a permutation of a stacked vector for all component correction currents $ \mathbf{I}^{\text{nl}}_c$.
Note that the voltage vector is augmented with the internal (auxiliary) voltage variables for the transformers (see Section~\ref{sec:transformers}). 

The variable current vector can now be stated as a function of the variable voltage vector, 
  \begin{IEEEeqnarray}{C} 
\mathbf{I}(\mathbf{U} ) = \mathbf{Y}\mathbf{U} \Rightarrow \mathbf{I}^\text{v}(\mathbf{U}^\text{v})
=
\mathbf{Y}^{\text{vf}} \mathbf{U}^\text{f} 
+
\mathbf{Y}^{\text{vv}} \mathbf{U}^\text{v}.
\label{eq_fix_var_decomp_ohms}  
\end{IEEEeqnarray}

\subsection{Guaranteeing Invertibility of $\mathbf{Y}^{\text{vv}}$}
Padding of impedance matrices with zeros to represent missing wires will lead to rank-deficient impedance matrices, so this is not allowed. 

For wye-to-delta transformers, it is important to pay attention the nature of this transformation matrix w.r.t. earth.
On the wye side, voltages are naturally defined w.r.t earth.
If the neutral is solidly grounded, the phase voltages are interpretable w.r.t. earth,
and if the neutral voltage is allowed to rise, a neutral voltage is interpretable w.r.t. to earth.
However, when traversing the transformer from the wye to the delta side,  the delta side only defines potential differences in its own framework, without referencing earth.
If any reference to earth is made downstream of the delta winding, the voltages will be naturally unique.
However, if no reference to earth is made, there remains a degree of freedom to resolve, also indicating that the admittance matrix is not invertible.
In real systems in the vicinity of earth, a reference to earth is unavoidable due to capacitive effects of conductors w.r.t earth.
Furthermore, line impedance matrices, when derived from Carson's equations, consider earth part of the circuit, so they imply the existence of such reference, even for three-wire networks where capacitance data has not been calculated. 

For these reasons, one can generally add a small shunt capacitance to components that do not have any capacitance defined in the input data.

\subsection{Current Injection Method Fixed Point Iteration Algorithm}
\begin{figure}[tbh]
  \centering
    \includegraphics[width=0.4\columnwidth]{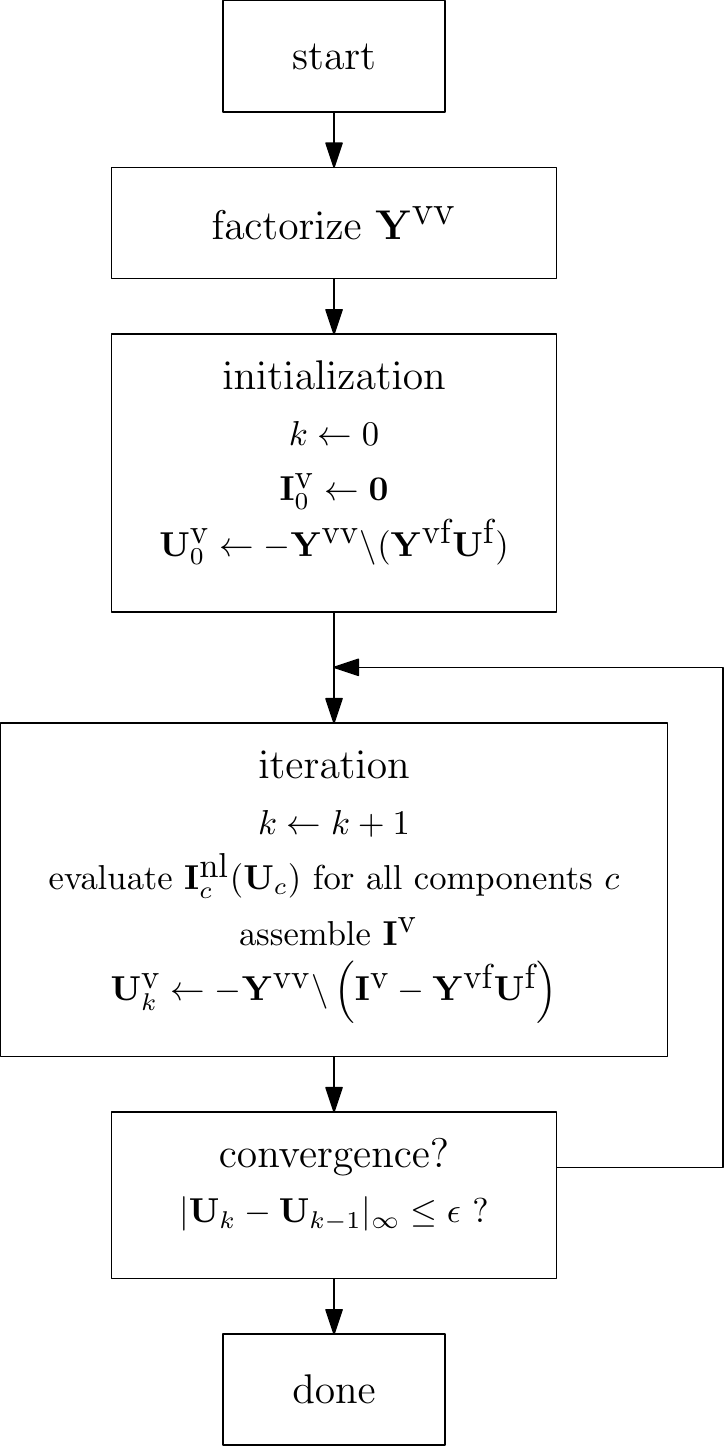}
  \caption{High-level flow of the current injection method fixed point iteration algorithm, performing only a single factorization. The inverses are implemented as matrix solves a.k.a. the `backslash' operation.}  \label{fig_fixedpointalgorihtm}
\end{figure}

Figure~\ref{fig_fixedpointalgorihtm} the high level flowchart of the \gls{CIM} algorithm with only a single factorization step. In the following we briefly explain how this flow is initialized, iterations, and the stopping criterion.

\subsubsection[Initialization]{Initialization $k=0$ }
We initialize the variable current vector at 0,
  \begin{IEEEeqnarray}{C} 
\mathbf{I}_0^\text{v}  \leftarrow  \mathbf{0}.
\end{IEEEeqnarray}
Therefore, we can write \eqref{eq_fix_var_decomp_ohms} as,
  \begin{IEEEeqnarray}{C} 
\mathbf{Y}^{\text{vv}} \mathbf{U}^{\text{v}}_0 = -\mathbf{Y}^{\text{vf}} \mathbf{U}^{\text{f}}
\end{IEEEeqnarray}
where $\mathbf{U}^{\text{f}}$ is the vector of known voltage and $\mathbf{U}^{\text{v}}_0$ is the initialization (iterate 0) of the vector of variable voltages.
We can now obtain the first estimates for the voltage by solving this equation,
  \begin{IEEEeqnarray}{C} 
 \mathbf{U}^{\text{v}}_0 \leftarrow  - (\mathbf{Y}^{\text{vv}})^{-1} \mathbf{Y}^{\text{vf}} \mathbf{U}^{\text{f}} .
\end{IEEEeqnarray}

Note that we do not actually invert the admittance matrix, we instead factorize it and perform multiple solves based on that single factorization.

\subsubsection[Iteration]{Iteration $k>0$}
We now solve \eqref{eq_fix_var_decomp_ohms} iteratively. First, we  evaluate all functions $\mathbf{I}^{\text{nl}}_c(\mathbf{U}^f,\mathbf{U}^v_{k-1})$ for the nonlinear components independently and then correctly place them in the permutated compensation current vector $\mathbf{I}^{\text{v}}$
  \begin{IEEEeqnarray}{C} 
\mathbf{Y}^{\text{vv}} \mathbf{U}^{\text{v}}_{k} = \mathbf{I}^{\text{v}}(\mathbf{U}^{\text{f}},\mathbf{U}^{\text{v}}_{k-1})-\mathbf{Y}^{\text{vf}} \mathbf{U}^{\text{f}}.
\end{IEEEeqnarray}
We can solve this equation for the next voltage estimate,
  \begin{IEEEeqnarray}{C} 
 \mathbf{U}^{\text{v}}_{k} \leftarrow  - (\mathbf{Y}^{\text{vv}})^{-1} \left(\mathbf{I}^{\text{v}}(\mathbf{U}^{\text{f}},\mathbf{U}^{\text{v}}_{k-1}) - \mathbf{Y}^{\text{vf}} \mathbf{U}^{\text{f}} \right).
\end{IEEEeqnarray}
Note that the matrix solve here is identical to that performed in the initialization, therefore we can re-use the stored factorization.

\subsubsection{Convergence Criterion}
When the voltage changes between consecutive iterations become sufficiently small, we have reached the fixed point. For this purpose we choose the infinity norm criterion,
  \begin{IEEEeqnarray}{C} 
\lVert \mathbf{U}_{k} - \mathbf{U}_{k-1} \rVert_{\infty} \leq \epsilon.
\end{IEEEeqnarray}
If the calculations are done in SI units, it is wise to do normalization to account for different voltage levels. 





\subsection{Computation Aspects}
We only need to factorize $\mathbf{Y}^{\text{vv}}$ once, during the initialization step. 
The factorization is then re-used every iteration.
It is possible to further minimize the number of iterations by re-building and re-factorizing the admittance matrix at each iteration. 
For instance, we re-evaluate the primitive admittance matrices for all constant-power loads in the new voltage phasor.  However, factorizing is computationally much more expensive than solving extra iteration steps. Therefore, updating the admittance matrix is best avoided, at the cost of an increase in the number of iterations.

\section{Concrete Element Types}
In this section we provided detailed derivation of the admittance matrix and compensation current for different concrete element types.

\subsection{Lines and Cables}
    Lines and cable are \emph{linear} power delivery elements, completely described by \eqref{eq_Y_lines} and \eqref{eq_Ic_lines} in \S\ref{sec_power_del_framework}.

\subsection{Switches}
Ideal (lossless) switches should really be removed from the circuit\footnote{e.g. through changing the input data or through Kron's reduction}. 
However, for validation purposes we replicate the series and shunt admittance matrices proposed in \textsc{OpenDSS} by treating them similar to the overhead lines and cables and assuming constant series and shunt matrices.

\subsection{Idealized and Lossy Transformers}
\label{sec:transformers}
We can decompose lossy n-winding transformers into a circuit composed of idealized single-phase transformers and series impedances and shunt admittances, using the methodology proposed in \cite{claeys2020tf}. 
We re-use the previously derived representations for lines and shunts, and add a component to represent idealized single-phase transformers.

\begin{figure}[tbh]
  \centering
    \includegraphics[width=0.40\columnwidth]{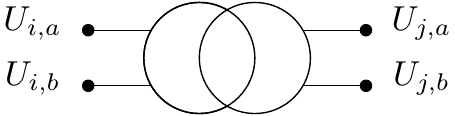}
  \caption{Idealized single-phase transformer.}  \label{fig_transformer}
\end{figure}

A single-phase idealized transformer is shown in Figure \ref{fig_transformer}.

We model the idealized transformer with an admittance matrix $\mathbf{Y}^{\text{tf}}$,
\begin{IEEEeqnarray}{C} 
\mathbf{I}_t =  \mathbf{Y}^{\text{tf}}_{t} \mathbf{U}_{t}.
 \end{IEEEeqnarray}
 
Table~\ref{table_transformers} provide the models for three transformer categories depending on the grounding configuration at either end. The model derivation is illustrated as follows:
\begin{itemize}
    \item With $r$ as the transformation ratio, we define the following admittance matrix of an idealized  transformer $\mathbf{Y}_{t}^{\text{tf}}$ for each type. To make the matrix invertible, we can augment the diagonal with a small positive scalar $\epsilon$.
    \item The $\mathbf{U}_{t}$ vector is the stack of the two from-side terminal voltages, the two to-side terminal voltages, and an auxiliary variable for the current, $I^{\text{aux}}_{t}$,
    \item The terminal currents are stacked and the last entry is set to zero.
    \item We solve the Ohm's law expressions and derive the expressions for the terminal currents, we derive the model summary with $ \epsilon \rightarrow 0$.
\end{itemize}

\begin{table}[tph]
\setlength{\tabcolsep}{3pt}
  \centering 
  \caption{Idealized single-phase transformer models.} 
  \label{table_transformers}
  \begin{tabular}{ m{2.2cm}  c  c  l }
    \toprule
     \\
     &
     $ \mathbf{Y}_{t}^{\text{tf}}$
     & 
     $\mathbf{I}_t, \mathbf{U}_{t} $
     & model summary
     \\
     \cmidrule(lr){2-2} \cmidrule(lr){3-3}  \cmidrule(lr){4-4} 
     \\
     No grounding at sending and receiving ends
     &
     $\begin{bmatrix}
        \cdot& \cdot& \cdot& \cdot&1/r\\
        \cdot&\cdot  &\cdot & \cdot& -1/r\\
        \cdot & \cdot&\cdot  & \cdot&-1 \\
        \cdot&\cdot & \cdot & \cdot&1 \\
        1/r& -1/r & -1 & 1 & \cdot\\
    \end{bmatrix} 
    + \epsilon \mathbf{I},$
     &
     $\begin{bmatrix}
        I_{tij,a} \\
        I_{tij,b} \\
        I_{tji,a} \\
        I_{tji,b} \\
        0
        \end{bmatrix},
        \begin{bmatrix}
        U_{i,a} \\
        U_{i,b} \\
        U_{j,a} \\
        U_{j,b} \\
        I^{\text{aux}}_{t}.
        \end{bmatrix}$
    &
    $\begin{cases}
        I_{tij,a} = (1/r) I^{\text{aux}}_{t},\\
        I_{tij,b} = - (1/r) I^{\text{aux}}_{t} ,\\
        I_{tji,a} = I^{\text{aux}}_{t} ,\\
        I_{tji,b} = - I^{\text{aux}}_{t}, \\
        U_{i,a} - U_{i,b} = r (U_{j,a} - U_{j,b} ).
    \end{cases}$
    \\
    \\
    Grounding at sending end
    &
    $\begin{bmatrix}
        \cdot&\cdot& \cdot&1/r\\
        \cdot  &\cdot  & \cdot&-1 \\
        \cdot  &\cdot  & \cdot&1 \\
        1/r&  -1 & 1& \cdot\\
    \end{bmatrix} 
    + \epsilon \mathbf{I},$
    &
    $\begin{bmatrix}
        I_{tij,a} \\
        I_{tji,a} \\
        I_{tji,b} \\
        0
        \end{bmatrix}, 
        \begin{bmatrix}
         U_{i,a} \\
        U_{j,a} \\
        U_{j,b} \\
        I^{\text{aux}}_{t}.
        \end{bmatrix}$
     &
     $\begin{cases}
        U_{i,a}  = r (U_{j,a} - U_{j,b} ), \\
        I_{tij,a} + I_{tji,a} /r = 0, \\
        I_{tji,a} + I_{tji,b} = 0
    \end{cases}$
    \\ 
    \\
    Grounding at sending and receiving ends
    &
    $\begin{bmatrix}
        \cdot& \cdot&1/r\\
        \cdot  & \cdot&-1 \\
        1/r&  -1 &  \cdot\\
        \end{bmatrix} 
        + \epsilon \mathbf{I},$
    &
    $\begin{bmatrix}
        I_{tij,a} \\
        I_{tji,a} \\
        0
        \end{bmatrix},
        \begin{bmatrix}
         U_{i,a} \\
         U_{j,a} \\
        I^{\text{aux}}_{t}
        \end{bmatrix}$
    &
    $\begin{cases}
        U_{i,a} = r U_{j,a},\\
        I_{tij,a} + I_{tji,a} /r =0
    \end{cases}$ 
    \\
    \\
     \bottomrule
  \end{tabular}
\end{table}

Note that the conservation of current is implied from these models. For instance, for the transformer with no grounding we have
\begin{IEEEeqnarray}{C} 
I_{tij,a} +I_{tij,b} = 0, \\
I_{tji,a} + I_{tji,b} =0, 
\end{IEEEeqnarray}
which satisfies
\begin{IEEEeqnarray}{C} 
I_{tij,a} = I_{tji,a}/r.
\end{IEEEeqnarray}
We also confirm this transformer is lossless since
\begin{multline}
S_{tij,a} + S_{tji,a} + S_{tij,b} + S_{tji,b} \\
= U_{i,a} I^*_{tij,a} +  U_{j,a} I^*_{tji,a} + U_{i,b} I^*_{tij,b} +  U_{j,b} I^*_{tji,b}  \\
= \underbrace{(U_{i,a} - U_{i,b})}_{r (U_{j,a} - U_{j,b} )} (\underbrace{I_{tij,a}}_{ I_{tji,a}/r })^* - (U_{j,a} -U_{j,b})I^*_{tji,a} = 0.
\end{multline}

\subsection{Loads and Generators - Wye Connected}
Wye-connected loads or generators are given 3 complex power setpoints, one between each phase and the neutral,
  \begin{IEEEeqnarray}{C} 
  \mathbf{S}_d^{\text{ref}} \leftarrow 
  \begin{bmatrix}
  S_{d,a}^{\text{ref}} \\
  S_{d,b}^{\text{ref}} \\
  S_{d,c}^{\text{ref}} \\
  \end{bmatrix}\in \mathbb{C}^{3\times 1}.
    \end{IEEEeqnarray}
  A reference voltage level is needed for all voltage-dependent models except for the  constant power model,
    \begin{IEEEeqnarray}{C} 
    \mathbf{U}_d^{\text{ref}}\leftarrow 
  \begin{bmatrix}
  U_{d,a}^{\text{ref}} \\
  U_{d,b}^{\text{ref}} \\
  U_{d,c}^{\text{ref}} \\
  \end{bmatrix}\in \mathbb{R}^{3\times 1}.
    \end{IEEEeqnarray}
    Using these parameters, we can define an equivalent admittance vector for each phase (no mutual) as
        \begin{IEEEeqnarray}{C} 
  \mathbf{Y}^{\text{ref}}_d  \leftarrow      (\mathbf{S}^{\text{ref}}_d)^* \oslash    (\mathbf{U}^{\text{ref}}_i)^{\circ2} 
  =\begin{bmatrix}
  (S^{\text{ref}}_{d,a})^* /(U^{\text{ref}}_{i,a})^2  \\
  (S^{\text{ref}}_{d,b})^* /(U^{\text{ref}}_{i,b})^2  \\
  (S^{\text{ref}}_{d,c})^* / (U^{\text{ref}}_{i,c})^2  \\
  \end{bmatrix}  \in \mathbb{C}^{3\times 1}.
  \end{IEEEeqnarray}
  
  The kron reduced model of the wye-connected load/generator has the admittance matrix of the form
  \begin{equation}
    \label{Yd_load_kron}
      \mathbf{Y}_d \leftarrow  
        \diagm(\mathbf{Y}^{\text{ref}}_d)   \in \mathbb{C}^{3\times 1},
  \end{equation}
  where as the model with explicit neutral has the admittance matrix,
  \begin{equation}
    \label{Yd_load_EN}
      \mathbf{Y}_d \leftarrow  
        \begin{bmatrix}
        \diagm(\mathbf{Y}^{\text{ref}}_d) &-\mathbf{Y}^{\text{ref}}_d \\
        -(\mathbf{Y}^{\text{ref}}_d)^\text{T} & \sum(\mathbf{Y}^{\text{ref}}_d)
        \end{bmatrix}   
        \in \mathbb{C}^{4\times 4}.
  \end{equation}
  
  Table~\ref{table_wye_loads_generators} lists compensation current vectors for kron-reduced and explicit-neutral variants of different load types (constant impedance, constant power, constant current, and exponential). For each load type, we need to derive the variable power and current vectors to calculate the compensation current vector.

  \begin{equation}
    \label{eq_Id_Sd}
      \begin{bmatrix}
      I_{d,a} \\
      I_{d,b} \\
      I_{d,c} \\
      \end{bmatrix}  
      \leftarrow 
      \begin{bmatrix}
      (S_{d,a}/ U_{i,a})^* \\
      (S_{d,b}/ U_{i,b})^*\\
      (S_{d,c}/ U_{i,c})^* \\
      \end{bmatrix}
  \end{equation}
  
  In models with explicit neutral we have to make sure to balance the current,
  \begin{IEEEeqnarray}{C} 
    \label{eq_Idn}
    I_{d,n} \leftarrow - \left(  I_{d,a} +   I_{d,b}  +   I_{d,c} \right).
  \end{IEEEeqnarray}

\begin{table}[tph]
\setlength{\tabcolsep}{3pt}
  \centering 
  \caption{Compensation currents for wye-connected loads and generators} 
  \label{table_wye_loads_generators}
  \begin{tabular}{ l  c c  c }
    \toprule 
     & Kron reduced && Explicit neutral      \\
     \cmidrule{2-2} \cmidrule{4-4} 
     \\
     Constant impedance
     & $\mathbf{I}^{\text{nl}}_d \leftarrow 0$
     && $\mathbf{I}^{\text{nl}}_d \leftarrow 0$
     \\ 
     \\
     \hline  & &
     \\
     Constant power
     &
     $\begin{aligned}[t]
      \begin{bmatrix}
      I_{d,a} \\
      I_{d,b} \\
      I_{d,c} \\
      \end{bmatrix}  
      &\leftarrow 
      \begin{bmatrix}
      ( S_{d,a}^{\text{ref}} / (U_{i,a} ) )^* \\
      ( S_{d,b}^{\text{ref}} / (U_{i,b} ) )^* \\
      ( S_{d,c}^{\text{ref}} / (U_{i,c} ) )^* \\
     \end{bmatrix}
     \\
     \mathbf{I}^{\text{nl}}_d 
     &\leftarrow  
     \mathbf{Y}_d\mathbf{U}_i  - \mathbf{I}_d   \in \mathbb{C}^{3\times 1}
     \end{aligned}$
     &&
     $\begin{aligned}[t]
     \begin{bmatrix}
      I_{d,a} \\
      I_{d,b} \\
      I_{d,c} \\
      \end{bmatrix}  
      &\leftarrow 
      \begin{bmatrix}
      ( S_{d,a}^{\text{ref}} / (U_{i,a} - U_{i,n} ) )^* \\
      ( S_{d,b}^{\text{ref}} / (U_{i,b} - U_{i,n} ) )^* \\
      ( S_{d,c}^{\text{ref}} / (U_{i,c} - U_{i,n} ) )^* \\
     \end{bmatrix}
     \\
     \mathbf{I}^{\text{nl}}_d 
     & \leftarrow  
     \mathbf{Y}_d\mathbf{U}_i  - \mathbf{I}_d   \in \mathbb{C}^{4\times 1}
     \end{aligned}$
    \\ 
    \\
    \hline & &
    \\
    Constant current
    &
    $\begin{aligned}[t]
     \begin{bmatrix}
      S_{d,a} \\
      S_{d,b} \\
      S_{d,c} \\
      \end{bmatrix}  
      &\leftarrow 
      \begin{bmatrix}
      S_{d,a}^{\text{ref}} \frac{ |U_{i,a}|}{  U_{d,a}^{\text{ref}} } \\
      S_{d,b}^{\text{ref}} \frac{ |U_{i,b}| }{  U_{d,b}^{\text{ref}} }\\
      S_{d,c}^{\text{ref}} \frac{ |U_{i,c}|}{  U_{d,c}^{\text{ref}} } \\
      \end{bmatrix}
     \\
     \mathbf{I}_d 
     &\leftarrow 
     \eqref{eq_Id_Sd}
     \\
     \mathbf{I}^{\text{nl}}_d
     &\leftarrow  
     \mathbf{Y}_d\mathbf{U}_i  - \mathbf{I}_d   \in \mathbb{C}^{3\times 1}
    \end{aligned}$
    &&
    $\begin{aligned}[t]
      \begin{bmatrix}
      S_{d,a} \\
      S_{d,b} \\
      S_{d,c} \\
      \end{bmatrix}  
      &\leftarrow 
      \begin{bmatrix}
      S_{d,a}^{\text{ref}} \frac{ |U_{i,a} - U_{i,n}|}{  U_{d,a}^{\text{ref}} } \\
      S_{d,b}^{\text{ref}} \frac{ |U_{i,b} - U_{i,n}| }{  U_{d,b}^{\text{ref}} }\\
      S_{d,c}^{\text{ref}} \frac{ |U_{i,c} - U_{i,n}|}{  U_{d,c}^{\text{ref}} } \\
      \end{bmatrix}
     \\
      \mathbf{I}_d 
      &\leftarrow 
      \eqref{eq_Id_Sd}, \eqref{eq_Idn}
     \\
     \mathbf{I}^{\text{nl}}_d 
     &\leftarrow  
     \mathbf{Y}_d\mathbf{U}_i  - \mathbf{I}_d   \in \mathbb{C}^{4\times 1}
    \end{aligned}$
   \\
   \\
   \hline & & 
   \\
   Exponential 
   &
   $\begin{aligned}[t]
    \begin{bmatrix}
      S_{d,a} \\
      S_{d,b} \\
      S_{d,c} \\
    \end{bmatrix}  
    &\leftarrow 
    \begin{bmatrix}
      P_{d,a}^{\text{ref}} \left(\frac{ |U_{i,a}|}{  U_{d,a}^{\text{ref}} } \right)^{x^P_{a}} \\
      P_{d,b}^{\text{ref}} \left(\frac{ |U_{i,b}| }{  U_{d,b}^{\text{ref}} }\right)^{x^P_{b}}\\
      P_{d,c}^{\text{ref}} \left(\frac{ |U_{i,c}|}{  U_{d,c}^{\text{ref}} } \right)^{x^P_{c}}\\
    \end{bmatrix} \\
    &+ j
    \begin{bmatrix}
      Q_{d,a}^{\text{ref}} \left(\frac{ |U_{i,a}|}{  U_{d,a}^{\text{ref}} } \right)^{x^Q_{a}} \\
      Q_{d,b}^{\text{ref}} \left(\frac{ |U_{i,b}|}{  U_{d,b}^{\text{ref}} }\right)^{x^Q_{b}}\\
      Q_{d,c}^{\text{ref}} \left(\frac{ |U_{i,c}|}{  U_{d,c}^{\text{ref}} } \right)^{x^Q_{c}}\\
    \end{bmatrix}
    \\
     \mathbf{I}_d 
     &\leftarrow 
     \eqref{eq_Id_Sd}
    \\
    \mathbf{I}^{\text{nl}}_d 
    &\leftarrow  
    \mathbf{Y}_d\mathbf{U}_i  - \mathbf{I}_d   \in \mathbb{C}^{3\times 1}
   \end{aligned}$
   &&
   $\begin{aligned}[t]
    \begin{bmatrix}
      S_{d,a} \\
      S_{d,b} \\
      S_{d,c} \\
    \end{bmatrix}  
    &\leftarrow 
    \begin{bmatrix}
      P_{d,a}^{\text{ref}} \left((\frac{ |U_{i,a} - U_{i,n}|}{  U_{d,a}^{\text{ref}} }\right)^{x^P_{a}} \\
      P_{d,b}^{\text{ref}} \left((\frac{ |U_{i,b} - U_{i,n}| }{  U_{d,b}^{\text{ref}} }\right)^{x^P_{b}}\\
      P_{d,c}^{\text{ref}} \left((\frac{ |U_{i,c} - U_{i,n}|}{  U_{d,c}^{\text{ref}} } \right)^{x^P_{c}}\\
    \end{bmatrix} \\
      &+ j
    \begin{bmatrix}
      Q_{d,a}^{\text{ref}} \left((\frac{ |U_{i,a} - U_{i,n}|}{  U_{d,a}^{\text{ref}} }\right)^{x^Q_{a}} \\
      Q_{d,b}^{\text{ref}} \left((\frac{ |U_{i,b} - U_{i,n}| }{  U_{d,b}^{\text{ref}} }\right)^{x^Q_{b}}\\
      Q_{d,c}^{\text{ref}} \left((\frac{ |U_{i,c} - U_{i,n}|}{  U_{d,c}^{\text{ref}} } \right)^{x^Q_{c}}\\
     \end{bmatrix}
    \\
      \mathbf{I}_d 
      &\leftarrow 
      \eqref{eq_Id_Sd}, \eqref{eq_Idn}
    \\ 
    \mathbf{I}^{\text{nl}}_d 
    &\leftarrow  
    \mathbf{Y}_d\mathbf{U}_i  - \mathbf{I}_d   \in \mathbb{C}^{4\times 1}
   \end{aligned}$
  \\
  \\
  \bottomrule
  \end{tabular}
\end{table}

\subsection{Loads and Generators - Delta Connected}
For delta-connected loads and generators, the power setpoints are interpreted as the product of the phase-to-phase voltages and the conjugate of the delta currents,
\begin{IEEEeqnarray}{C} 
  \mathbf{S}_d^{\Delta,\text{ref}} =
  \begin{bmatrix}
  S_{d,ab}^{\text{ref}} \\
  S_{d,bc}^{\text{ref}} \\
  S_{d,ca}^{\text{ref}} \\
  \end{bmatrix}
  = \mathbf{U}^{\Delta}_i \circ \left (\mathbf{I}^{\Delta}_d \right)^*
  \in \mathbb{C}^{3\times 1} .
\end{IEEEeqnarray}
We define a reference value for the phase-to-phase voltage magnitude. For instance, in a $230$\,V nominal phase voltage network, this would be a length-3 vector filled with the value $230 \sqrt{3}$\,V,
\begin{IEEEeqnarray}{C} 
  \mathbf{U}_d^{\Delta,\text{ref}}\leftarrow 
  \begin{bmatrix}
  U_{d,ab}^{\text{ref}} \\
  U_{d,bc}^{\text{ref}} \\
  U_{d,ca}^{\text{ref}} \\
  \end{bmatrix}\in \mathbb{R}^{3\times 1}.
\end{IEEEeqnarray}
    We then derive the vector with the equivalent shunt admittance,
\begin{IEEEeqnarray}{C} 
  \mathbf{Y}^{\Delta,\text{ref}}_d  \leftarrow    \left(  \mathbf{S}^{\Delta,\text{ref}}_d \right)^* \oslash    \left(\mathbf{U}^{\Delta,\text{ref}}_i  \right)^{\circ2} 
  =\begin{bmatrix}
  (S^{\text{ref}}_{d,ab})^* /(U^{\text{ref}}_{i,ab})^2  \\
  (S^{\text{ref}}_{d,bc})^* /(U^{\text{ref}}_{i,bc})^2  \\
  (S^{\text{ref}}_{d,ca})^* / (U^{\text{ref}}_{i,ca})^2  \\
  \end{bmatrix}  \in \mathbb{C}^{3\times 1},
\end{IEEEeqnarray}
where superscript $\circ2$ indicates the element-wise square and $\oslash$ the element-wise division.
The line voltage vector is a linear transformation of the phase voltage vector, with transformation matrix $\mathbf{M}$,
  \begin{IEEEeqnarray}{C} 
\mathbf{U}^{\Delta}_i = 
\begin{bmatrix}
U_{i,ab} \\
U_{i,bc} \\
U_{i,ca} \\
\end{bmatrix} 
=
 \begin{bmatrix}
\phantom{-}1 & -1 & \phantom{-}0 \\
\phantom{-}0 & \phantom{-}1 & -1  \\
-1 & \phantom{-}0 & \phantom{-}1   \\
\end{bmatrix}
\begin{bmatrix}
U_{i,a} \\
U_{i,b} \\
U_{i,c} \\
\end{bmatrix} 
=
\mathbf{M} \mathbf{U}_i . \label{eq_voltage_delta_def}
\end{IEEEeqnarray}
Such constant-admittance delta loads satisfy,
\begin{IEEEeqnarray}{C} 
  \mathbf{I}^{\Delta}_d = \mathbf{Y}^{\Delta,\text{ref}}_d \circ \mathbf{U}^{\Delta}_i. \label{eq_delta_load_def}
\end{IEEEeqnarray}

Similarly, the phase current vector $ \mathbf{I}_c $ is a linear transformation of the line delta current vector $ \mathbf{I}^{\Delta}_d$, but with transformation matrix $\mathbf{M}^{\text{T}}  $,
\begin{IEEEeqnarray}{C} 
   \mathbf{I}_c = 
\begin{bmatrix}
I_{c,a} \\
I_{c,b} \\
I_{c,c} \\
\end{bmatrix} = 
 \begin{bmatrix}
\phantom{-}1 & \phantom{-}0 & -1 \\
-1 & \phantom{-}1 & \phantom{-}0 \\
\phantom{-}0 & -1 & \phantom{-}1   \\
\end{bmatrix}
   \begin{bmatrix}
  I_{d,ab} \\
  I_{d,bc} \\
  I_{d,ca} \\
  \end{bmatrix}  
  = 
\mathbf{M}^{\text{T}}    \mathbf{I}^{\Delta}_d . \label{eq_bus_phase_delta_current}
\end{IEEEeqnarray}
We now substitute \eqref{eq_delta_load_def} into \eqref{eq_bus_phase_delta_current}, and use \eqref{eq_voltage_delta_def} to derive the appropriate admittance matrix for delta-connected load and generators, 
    \begin{IEEEeqnarray}{C} 
  \mathbf{Y}_d \leftarrow  
\mathbf{M}^{\text{T}}\diagm(\mathbf{Y}^{\Delta,\text{ref}}_d)  \mathbf{M}.
\end{IEEEeqnarray}

Table~\ref{table_delta_loads_generators} lists the compensation current vectors for different load types. Again, the variable current vector should be derived to compute the compensation current vector.

\begin{table}[tph]
\setlength{\tabcolsep}{3pt}
  \centering 
  \caption{Compensation currents for delta-connected loads and generators} 
  \label{table_delta_loads_generators}
  \begin{tabular}{ l  l }
    \toprule
     \\
     Constant impedance
     & $\mathbf{I}^{\text{nl}}_d \leftarrow 0$
     \\ 
     \\
     \hline  &
     \\
     Constant power
     &
     $\begin{aligned}[t]
        & \mathbf{I}^{\Delta}_d = 
        \begin{bmatrix}
          I_{d,ab} \\
          I_{d,bc} \\
          I_{d,ca} \\
        \end{bmatrix}  
        \leftarrow 
        \begin{bmatrix}
        ( S_{d,ab}^{\text{ref}} / (U_{i,a} - U_{i,b} ) )^* \\
        ( S_{d,bc}^{\text{ref}} / (U_{i,b} - U_{i,c} ) )^* \\
        ( S_{d,ca}^{\text{ref}} / (U_{i,c} - U_{i,a}  ) )^* \\
        \end{bmatrix}
      \\
       &  \mathbf{I}^{\text{nl}}_d \leftarrow  \mathbf{Y}_d\mathbf{U}_i  - \mathbf{M}^{\text{T}}    \mathbf{I}^{\Delta}_d 
     \end{aligned}$
     \\ 
     \\
     \hline  &
     \\
     Constant current
     &
     $\begin{aligned}[t]
        &\qquad \ 
        \begin{bmatrix}
          S_{d,a} \\
          S_{d,b} \\
          S_{d,c} \\
        \end{bmatrix}  
        \leftarrow 
          \begin{bmatrix}
          S_{d,a}^{\text{ref}} \frac{ |U_{i,a} - U_{i,b}|}{  U_{d,a}^{\text{ref}} } \\
          S_{d,b}^{\text{ref}} \frac{ |U_{i,b} - U_{i,c}|}{  U_{d,b}^{\text{ref}} }\\
          S_{d,c}^{\text{ref}} \frac{ |U_{i,c} - U_{i,a}|}{  U_{d,c}^{\text{ref}} } \\
        \end{bmatrix}
        \\
        & \mathbf{I}^{\Delta}_d = 
        \begin{bmatrix}
          I_{d,ab} \\
          I_{d,bc} \\
          I_{d,ca} \\
        \end{bmatrix}  
        \leftarrow 
        \begin{bmatrix}
        ( S_{d,ab} / (U_{i,a} - U_{i,b} ) )^* \\
        ( S_{d,bc} / (U_{i,b} - U_{i,c} ) )^* \\
        ( S_{d,ca} / (U_{i,c} - U_{i,a}  ) )^* \\
        \end{bmatrix}
      \\
       &  \mathbf{I}^{\text{nl}}_d \leftarrow  \mathbf{Y}_d\mathbf{U}_i  - \mathbf{M}^{\text{T}}    \mathbf{I}^{\Delta}_d 
     \end{aligned}$
    \\
    \\
    \hline  &
    \\
    Exponential
     &
     $\begin{aligned}[t]
        &\qquad \ 
        \begin{bmatrix}
          S_{d,a} \\
          S_{d,b} \\
          S_{d,c} \\
        \end{bmatrix}  
        \leftarrow 
        \begin{bmatrix}
          P_{d,a}^{\text{ref}} \left(\frac{ |U_{i,a}-U_{i,b}|}{  U_{d,a}^{\text{ref}} } \right)^{x^P_{a}} \\
          P_{d,b}^{\text{ref}} \left(\frac{ |U_{i,b}-U_{i,c}| }{  U_{d,b}^{\text{ref}} }\right)^{x^P_{b}}\\
          P_{d,c}^{\text{ref}} \left(\frac{ |U_{i,c}-U_{i,a}|}{  U_{d,c}^{\text{ref}} } \right)^{x^P_{c}}\\
        \end{bmatrix} \\
         & \qquad \qquad \qquad + j
        \begin{bmatrix}
          Q_{d,a}^{\text{ref}} \left(\frac{ |U_{i,a}-U_{i,b}|}{  U_{d,a}^{\text{ref}} } \right)^{x^Q_{a}} \\
          Q_{d,b}^{\text{ref}} \left(\frac{ |U_{i,b}-U_{i,c}|}{  U_{d,b}^{\text{ref}} }\right)^{x^Q_{b}}\\
          Q_{d,c}^{\text{ref}} \left(\frac{ |U_{i,c}-U_{i,a}|}{  U_{d,c}^{\text{ref}} } \right)^{x^Q_{c}}\\
        \end{bmatrix}
        \\
        & \mathbf{I}^{\Delta}_d = 
        \begin{bmatrix}
          I_{d,ab} \\
          I_{d,bc} \\
          I_{d,ca} \\
        \end{bmatrix}  
        \leftarrow 
        \begin{bmatrix}
        ( S_{d,ab} / (U_{i,a} - U_{i,b} ) )^* \\
        ( S_{d,bc} / (U_{i,b} - U_{i,c} ) )^* \\
        ( S_{d,ca} / (U_{i,c} - U_{i,a}  ) )^* \\
        \end{bmatrix}
      \\
       &  \mathbf{I}^{\text{nl}}_d \leftarrow  \mathbf{Y}_d\mathbf{U}_i  - \mathbf{M}^{\text{T}}    \mathbf{I}^{\Delta}_d 
     \end{aligned}$
    \\
    \\
    \bottomrule
  \end{tabular}
\end{table}




\section{Implementation and Validation}
\subsection{Differences w.r.t. \textsc{OpenDSS}}

This report provides a summary of the power flow solver implemented in \textsc{PowerModelsDistribution.jl} (PMD),  inspired by the existing open-source solver included in \textsc{OpenDSS}. However, there are still some differences between the two implementations:
\begin{itemize}
    \item The \gls{PMD} implementation is valid for up to four-wire networks, but not more, whereas \textsc{OpenDSS} can solve cases with more wires. Supporting n-wire should be relatively straight-forward though, as most of the code is agnostic, so the key challenge is the  validation of such an extension. 
    \item \gls{PMD} computes power flow results in per unit values whereas \textsc{OpenDSS} does not\footnote{For the \emph{optimization} models in \gls{PMD}, the per-unit system is overall beneficial in the authors' experience. Simultaneously, the authors agree with Roger Dugan that the per unit system is not beneficial for simulation problems. Nevertheless, due to code design it was faster to prototype \gls{CIM}s native solver starting from the internal per-unit representation.}. The per unit system interacts with the addition of small diagonal values. E.g. OpenDSS adds the equivalent of 10 kvar capacitive at 345 kV\footnote{Available: \url{https://github.com/tshort/OpenDSS/blob/master/Source/PDElements/Line.pas}}, which is easy to do in SI units in the matrix manipulation code context, but requires additional context if it is performed in per-unit. 
    \item The load model relaxation feature in \textsc{OpenDSS} (changing the load model to constant impedance outside of the \texttt{vminpu, vmaxpu} range) has not been implemented in \gls{CIM}. The lack of this relaxation makes convergence more challenging. 
    \item For maximum compatibility, \gls{CIM} uses Julia's default matrix factorization routines whereas \textsc{OpenDSS} uses the KLU library. KLU and other specialized matrix solvers could be accessed from Julia through \textsc{LinearSolve.jl}\footnote{Available: \url{https://github.com/SciML/LinearSolve.jl}}.
    \item The transformer decomposition described above \cite{claeys2020tf} and implemented in \gls{CIM} is slightly different from how OpenDSS models transformers (additional internal nodes).
    \item \textsc{OpenDSS} adds small values (ppm) to the diagonal of the admittance matrix for various elements and \gls{CIM} does not replicate all those cases.
\end{itemize}

\subsection{Validation}
The \gls{CIM} power flow results are validated against \textsc{OpenDSS} results (calculated using \textsc{OpenDSSDirect.jl}\footnote{Available: \url{https://github.com/dss-extensions/OpenDSSDirect.jl}}  0.8.1) for which Table~\ref{table_validation_results} lists the maximum voltage error in per unit between the \gls{CIM} (tolerance 1E-8) and OpenDSS (tolerance 1E-10) as
\begin{equation}
    \mathbf{U}_{\max}^{pu} = \max_{i \in \mathcal{I}, p \in \mathcal{P}} |\mathbf{U}_{i,p}^{\text{CIM}} - \mathbf{U}_{i,p}^{\text{ODSS}} |
\end{equation}
where $\mathbf{U}_{i,p}^{\text{CIM}}$ and $\mathbf{U}_{i,p}^{\text{ODSS}}$ are respectively bus terminal voltages calculated by \gls{CIM} and \textsc{OpenDSSDirect.jl}.
The first batch of test cases are unit tests of the \gls{CIM} package, and the second batch are larger networks.
We aim to update this table when changes are made, or more cases are validated. 

\paragraph{Note:} It should also be mentioned that the implementation still has room for improvement as some cases still result in singularity errors. 
For instance the \texttt{Egrid - SantaFe/urban-suburban/uhs0 1247 uhs0 1247 -udt4077} network.

\begin{landscape}
\thispagestyle{mylandscape} 
\begin{table}[htbp]
    \centering
    {\footnotesize
    \begin{tabular}{l c c c c c c}
         \toprule
         & \multicolumn{4}{c}{\textsc{PowerModelsDistribution}} & \multicolumn{2}{c}{\textsc{OpenDSS}}  \\
          \cmidrule(lr){2-5}   \cmidrule(lr){6-7} \\
         Network & $\mathbf{U}_{\max}^{pu}$ & time total (s) & iterations & $\mathbf{Y}^{\text{vv}}$ size & iterations & time (s) \\
         &  & (build/solve/post) & & &\\
         \hline 
         1-phase EN*  -- wye generator & 1.1319E-9 & 0.00106 & 8 & 7 $\times$ 7 & 6 & 0.0006 \\
         1-phase EN  -- delta generator & 3.9545E-10 & 0.00094 & 6 & 7 $\times$ 7 & 7 & 0.0007 \\
         3-phase EN  -- wye generator & 1.1375E-9 & 0.00133 & 9 & 7 $\times$ 7 & 9 & 0.0007 \\
         3-phase EN  -- delta generator & 1.3978E-9 & 0.00101 & 7 & 7 $\times$ 7 & 7 & 0.0007 \\
         1-phase EN  -- wye constant power load & 8.9164E-8 & 0.00184 & 31 & 7 $\times$ 7 & 15 & 0.0011\\
         1-phase EN  -- delta constant power load & 1.2511E-9 & 0.00095 & 8 & 7 $\times$ 7 & 7 & 0.0007 \\
         3-phase EN  -- wye constant impedance load & 1.2831E-9 & 0.00097 & 8 & 7 $\times$ 7  & 7 & 0.0006 \\
         3-phase EN  -- wye constant current load & 4.8127E-10 & 0.00114 & 10 & 7 $\times$ 7  & 8 & 0.0007 \\
         3-phase EN  -- wye constant power load & 1.7351E-9 & 0.00110 & 11 & 7 $\times$ 7  & 9 & 0.0007 \\
         3-phase EN  -- delta constant impedance load & 3.3809E-10 & 0.00101 & 9 & 7 $\times$ 7 & 7 & 0.0007 \\
         3-phase EN  -- delta constant current load & 3.6498E-10 & 0.00123 & 9 & 7 $\times$ 7 & 7 & 0.0007 \\
         3-phase EN  -- delta constant power load & 5.2944E-10 & 0.00127 & 9 & 7 $\times$ 7 & 7 & 0.0007 \\
         3-phase EN  -- switch & 4.7188E-8 & 0.00128 & 8 & 19 $\times$ 19 & 7 & 0.0008 \\
         3-phase  -- switch & 3.3801E-10 & 0.00117 & 6 & 15 $\times$ 15 & 5 & 0.0007 \\
         1-phase  -- switch & 1.4997E-8 & 0.00073 & 6 & 6 $\times$ 6 & 5 & 0.0007 \\
         3-phase  -- 3-winding (yyy) transformer & 7.2649E-7 & 0.00142 & 7 & 17 $\times$ 17 & 4 & 0.0024 \\
         3-phase  -- 3-winding (dyy) transformer & 7.4705E-7 & 0.00233 & 15 & 39 $\times$ 39 & 8 & 0.0008 \\
         3-phase  -- 2-winding (dy) transformer & 5.8773E-9 & 0.00128 & 3 & 27 $\times$ 27 & 2 & 0.0007 \\
         3-phase  -- 2-winding (yy) transformer & 5.9322E-9 & 0.00134 & 3 & 27 $\times$ 27 & 2 & 0.0007 \\
         3-phase EN  -- 2-winding (dy) transformer & 6.6057E-8 & 0.00239 & 10 & 33 $\times$ 33 & 8 & 0.0009 \\
         3-phase EN  -- 2-winding (yy) transformer & 7.4026E-8 & 0.00181 & 10 & 33 $\times$ 33 & 8 & 0.0010\\
         \hline
         IEEE 13 3-phase  &  3.7651E-6 & 0.00726 & 15 & 89 $\times$ 89 & 9 & 0.0022 \\
          & & (0.0037/0.0026/0.0000) & & & \\
         IEEE 34 3-phase  &  6.8018E-8 & 0.03098 & max & 210 $\times$ 210 & 5 & 0.0040 \\
          & & (0.0116/0.0162/0.0032) & & &   \\
         IEEE 123 3-phase  & 4.0450E-8 & 0.03162 & 28 & 343 $\times$ 343 & 6 & 0.0068 \\
          & & (0.0116/0.0167/0.0033) & & & \\
         Egrid - GreensBoro/Industrial  & 0.0018  & 26.9685 & max & 46697 $\times$ 46697 & 6 & 0.4497 \\
          & & (21.0289/5.6004/0.3392) & & \\
         Egrid - SantaFe/urban-suburban/uhs0 1247/ &   0.00012 & 26.4347 & max & 13047 $\times$ 13047 & 6 & 0.1696\\ 
         \qquad uhs0 1247 -udt4776 & & (20.5699/5.45653/0.4083) & & \\
         \bottomrule 
    \end{tabular}
    }
    \caption{\textsc{PowerModelsDistribution} unbalanced power flow (with tolerance 1E-8) validation results w.r.t to \textsc{OpenDSS} (with tolerance 1E-10). \\ $^*$ EN: with explicit neutral.}
    \label{table_validation_results}
\end{table}
\end{landscape}

\subsection{Future work}
The implementation can still be improved, and we invite collaborators to work on these issues:
\begin{itemize}
    \item implement the load model relaxation as used in OpenDSS, allowing faster convergence on more test cases;
    \item post-processing of congestions (tagging of occurrences of under/over voltage/ current, deriving of reliability metrics);
    \item testing limits of scalability of the implementation, e.g. across a variety of linear solvers;
    \item continue validation efforts on more distribution network test cases;
    \item alternative models for switches, e.g. Kron-reducing buses as part of the algorithm;
    \item implementation without per unit conversion, without transformer decomposition;
    \item volt-var/volt-watt control for inverters, including support for arbitrary control laws and multiple voltage inputs\footnote{E.g. OpenDSS only supports voltage-to-ground inputs for the volt-var/watt control laws};
    \item multi-time-step support that re-uses the previous solution and factorized nodal admittance matrix where possible;
    \item uncertainty extensions, e.g. polynomial chaos or Monte-Carlo simulation.
\end{itemize}

As a final note, there is a need for more flexible distribution network data models that can capture unbalanced power flow as well as unbalanced OPF and state estimation data sets with well-defined semantics and libraries that check for issues and inconsistencies in such data sets. 
For instance, in OpenDSS one can mix Kron-reduced-neutral and explicit-neutral representations in a single case study, which makes data checks difficult. 
Furthermore, in the author's opinions, load model relaxation and Kron's reduction should be established as settings of the power flow algorithm, not features of a data set describing a network.

\bibliographystyle{IEEEtran}
\bibliography{library,library2}

\end{document}